\documentclass[10pt,twoside]{classe-ma_e}

\usepackage{amscd}      
\usepackage{amssymb}     
\usepackage[intlimits]{amsmath}     
\usepackage{xypic}      
\LaTeXdiagrams          
\usepackage[all,v2]{xy}      
\xyoption{2cell}
\UseAllTwocells
\xyoption{frame}
\CompileMatrices

            {\nolinebreak $\Box$ \end{trivlist}}

\newcommand{\noprint}[1]{}

\newcommand{\toto}{\rightrightarrows}

\newcommand{\upst}{^{\ast}}

\newcommand{\com}{^{\scriptscriptstyle\bullet}}
\newcommand{\lcom}{_{\scriptscriptstyle\bullet}}

\newcommand{\XX}{{\mathfrak X}}

\newcommand{\PP}{{\mathfrak P}}

\newcommand{\GG}{{\mathfrak G}}

\newcommand{\zz}{{\mathbb Z}}

\newcommand{\rr}{{\mathbb R}}

\newcommand{\del}{\partial}

\newcommand{\id}{\mathop{\rm id}\nolimits}

\newcommand{\ldiag}[1]%
       {\makebox[0cm]{${\scriptstyle#1}\downarrow\phantom{\scriptstyle#1}$}}
\newcommand{\ldiagup}[1]%
       {\makebox[0cm]{${\scriptstyle#1}\uparrow\phantom{\scriptstyle#1}$}}
\newcommand{\rdiag}[1]%
       {\makebox[0cm]{$\phantom{\scriptstyle#1}\downarrow{\scriptstyle#1}$}}
\newcommand{\sediagr}[1]%
       {\makebox[0cm]{$\phantom{\scriptstyle#1}\searrow{\scriptstyle#1}$}}
\newcommand{\nediagr}[1]%
       {\makebox[0cm]{$\phantom{\scriptstyle#1}\nearrow{\scriptstyle#1}$}}
\newcommand{\rdiagup}[1]%
       {\makebox[0cm]{$\phantom{\scriptstyle#1}\uparrow{\scriptstyle#1}$}}
\newcommand{\swdiag}[1]%
       {\makebox[0cm]{$\phantom{\scriptstyle#1}\swarrow{\scriptstyle#1}$}}
\newcommand{\sediag}[1]%
       {\makebox[0cm]{${\scriptstyle#1}\searrow\phantom{\scriptstyle#1}$}}
\newcommand{\nediag}[1]%
       {\makebox[0cm]{${\scriptstyle#1}\nearrow\phantom{\scriptstyle#1}$}}

\newcommand{\longiso}{\stackrel{\textstyle\sim}{\longrightarrow}}

\newcommand{\doublearrowstack}[2]%
                      {{{{\scriptstyle#1}\atop{\textstyle\longrightarrow}}\atop{{\textstyle\longrightarrow}\atop{\scriptstyle#2}}}}
\newcommand{\rightleftarrowstack}[2]%
                      {{{{\scriptstyle#1}\atop{\textstyle\longrightarrow}}\atop{{\textstyle\longleftarrow}\atop{\scriptstyle#2}}}}
\newcommand{\leftrightarrowstack}[2]%
                      {{{{\scriptstyle#1}\atop{\textstyle\longleftarrow}}\atop{{\textstyle\longrightarrow}\atop{\scriptstyle#2}}}}

\newcommand{\overtoparrow}%
{\makebox[0cm]{\beginpicture
\setcoordinatesystem units <.8cm,.4cm> point at 0 0
\setplotarea x from -3 to 3, y from 0 to 1
\setquadratic
\plot -3 0 0 1 3 0 /
\put{\vector(3,-1){0}}[Bl] at 3 0
\endpicture}}

\newcommand{\underbottomarrow}%
{\makebox[0cm]{\beginpicture
\setcoordinatesystem units <.8cm,.4cm> point at 0 0
\setplotarea x from -3 to 3, y from 0 to 1
\setquadratic
\plot -3 1 0 0 3 1 /
\put{\vector(3,1){0}}[Bl] at 3 1
\endpicture}}

\newcommand{\ses}[5]%
{0\longrightarrow#1\stackrel{#2}{ \longrightarrow}#3\stackrel{#4}{
\longrightarrow}#5\longrightarrow0}

\newcommand{\dt}[6]%
{#1\stackrel{#2}{longrightarrow}#3 \stackrel{#4}{\longrightarrow}#5
\stackrel{#6}{\longrightarrow} #1[1]}  
 
\newcommand{\cat}[1]%
{(\mbox{\rm #1})}

\newcommand{\gm}{\Gamma}

\def\gpd{\,\lower1pt\hbox{$\longrightarrow$}\hskip-.24in\raise2pt
             \hbox{$\longrightarrow$}\,}

\newtheorem{them}{Theorem}[section]

\newtheorem{prop}[them]{Proposition}

\newtheorem{defn}[them]{Definition\rm}
\newtheorem{numrmk}{\it Remark\/}

\ComParit{S0764-4442}
\PIT{FLA}
\PXMA{????}
\Add{?}
\Volume{332}
\Year{2001}
\FirstPage{1}
\LastPage{??}
\AuteurCourant{Running authors}
\TitreCourant{Running title} 
\Journal
\Rubrique{Rubrique}{Heading}
\SousRubrique{Sous-rubrique}{Sub-Heading}
\PresentePar{First name}{NAME}
\Recu{jour mois ann\'ee}{apr\`es r\'evision}{jour mois ann\'ee}
\begin{document}
\title{%
$S^1$-Bundles and Gerbes over Differentiable Stacks
}
\author{%
Kai Behrend~$^{\text{a}}$,\ \
Ping Xu~$^{\text{b}}$
}
\address{%
\begin{itemize}\labelsep=2mm\leftskip=-5mm
\item[$^{\text{a}}$]
University of British Columbia \\
E-mail: behrend@math.ubc.ca
\item[$^{\text{b}}$]
Pennsylvania State University \\
E-mail: ping@math.psu.edu
\end{itemize}
}
\maketitle
\thispagestyle{empty}
\begin{Abstract}{%
We study $S^1$-bundles and $S^1$-gerbes over
differentiable stacks in terms of Lie groupoids, and construct
Chern classes and Dixmier-Douady classes in terms of analogues of
connections and curvature. 
}\end{Abstract}
\begin{Ftitle}{%
$S^1$-Fibr\'es et Gerbes sur des Champs Diff\'erentiables
}\end{Ftitle}
\begin{Resume}{%
On \'etudie les $S^1$-fibr\'es et les $S^1$-gerbes sur des champs diff\'erentiables
en termes de groupo\"\i des de Lie et construit les classes de Chern et
Dixmier-Douady en
termes d'analogues aux connexions et courbure. 
}\end{Resume}

\AFv

Soit $\XX$ un champ diff\'erentiable et $\PP\to\XX$ un $S^1$-fibr\'e
sur $\XX$. Soit $\Gamma\toto M$ une pr\'esentation par un groupo\"\i de de Lie pour 
$\XX$. Alors $\PP$ induit un $S^1$-fibr\'e $P$ sur $M$ sur lequel
agit $\Gamma\toto M$. On r\'ealise la classe de Chern de $\PP$ en
termes de don\'ees de type connexion sur $P$ et	prouve l'existance
des pr\'equantifications. Plus pr\'ecis\'ement, Soit $\theta\in \Omega^1(P)$ une 
pseudo-connexion, et $\omega+\Omega\in Z^2_{DR}(\Gamma\lcom)$ sa
pseudo-courbure. 

\begin{them}
La classe $[\omega+\Omega]\in H^2_{DR}(\Gamma\lcom)$ est
ind\'ependante du choix de la pseudo-connexion $\theta$ et 
correspond \`a la classe de Chern de $P$.
R\'eciproquement, soit $\omega+\Omega\in C^2_{DR} (\Gamma\lcom)$ un 
2-cocycle entier. Alors il existe un $S^1$-fibr\'e $P$
sur $\Gamma\toto M$ et une pseudo-connexion $\theta\in \Omega^1(P)$
ayant $\omega+\Omega$ pour pseudo-courbure.
De plus, l'ensemble des classes d'isomorphisme de tous ces $(P,\theta)$ est un
 $H^1 (\Gamma\lcom , \rr/\zz )$-ensemble.
\end{them}

Si $\GG$ est une $S^1$-gerbe sur $\XX$, et $R\toto M$ une pr\'esentation 
du champ diff\'erentiable $\GG$ et soit $\Gamma\toto M$ le groupo\"\i de de Lie
d\'efini par la pr\'esentation induite $M\to \XX$ de $\XX$.
Alors $R$ est une $S^1$-extention centrale du groupo\"\i de de Lie $\Gamma\toto M$.
Ainsi les $S^1$-extensions centrales de $\Gamma\toto M$ sont exactement les 
$S^1$-gerbes sur $\XX$, donn\'ees d'une trivialisation sur $M$.
A nouveau, on peut r\'ealiser les classes caract\'eristiques de la gerbe (que nous appelons
classes de Dixmier-Douady) en termes de donn\'ees de type connexion et prouver
l'existances de pr\'equantifications . Plus pr\'ecis\'ement, soit 
$\theta+B \in C^2_{DR} (R)$ une pseudo-connexion sur
$R$, et $\theta+\omega+\Omega \in Z^3_{DR}(\Gamma\lcom)$ sa 
pseudo-courbure. 

\begin{them}
La classe $[\eta+\omega +\Omega ]\in H^3_{DR}(\Gamma\lcom)$ est ind\'ependante du
choix de la pseudo-connexion $\theta +B$ sur $R$ et
correspond \`a la classe de Dixmier-Douady de $R$. 
R\'eciproquement, pour tout 3-cocycle $\eta +\omega+\Omega\in
Z^3_{DR}(\Gamma\lcom)$
tel que $[\eta+\omega+\Omega]$ est une classe enti\`ere et
$\Omega$ est exact,
il existe une extension centrale $R\toto M$ du
groupo\"\i de $\gm \toto M$, et une pseudo-connexion $\theta+B \in C^2_{DR} (R)$
sur $R$ telle que $\eta+\omega+\Omega$ soit la pseudo-courbure. 
Les paires $(R, \theta, B)$ forment, \`a un isomorphisme pr\`es, un ensemble simplement transitif 
sous le groupe des extensions centrales plates.

\end{them}
Dans le cas s-connexe, on obtient un construction explicite de
l'extension centrale avec pseudo-connexion. Cela donne
\'egalement un crit\`ere pour qu'une classe dans $H^3_{DR}(\Gamma\lcom)$ soit
enti\`ere. Ce th\'eor\`eme g\'en\'eralise le r\'esultat de~\cite{WX}. 

\begin{them}
Soit $\gm \toto M$ un groupo\"\i de de Lie $s$-connexe, et $\eta +\omega
\in C^3_{DR} (\gm\lcom )$ un 3-cocycle, o\`u $\eta \in \Omega^1(\gm_2
)$ and $\omega \in \Omega^2 (\gm )$. Supposons que $\omega $ repr\'esente
une classe de cohomologie enti\`ere dans $H^2_{DR}(\Gamma)$, de telle sorte qu'il
existe un $S^1$-fibr\'e $\pi:R\to \gm$ avec une connexion $\theta \in
\Omega^1 (R)$, dont la courbure est $\omega$.
Supposons que $\epsilon\upst R$, dot\'e d'une connexion plate
$\epsilon\upst\theta+\pi\upst \epsilon_2\upst\eta$ soit sans holonomie. (Ici
$\epsilon:M\to\Gamma$ et $\epsilon_2:M\to \Gamma_2$ sont les morphisme d'identit\'e
respectifs.)
Alors $R\toto M$ admet de façon naturelle une structure de groupo\"\i de,
telle que $R$ soit une extension $S^1$-centrale de $\Gamma\toto M$
et $\eta +\omega$ la pseudo-courbure de $\theta$.
\end{them} 
 
Puisque les extensions centrales de groupo\"\i des d\'ecrivent les gerbes sur $\XX$
avec des trivialisations donn\'ees sur $M$, on peut seulement d\'ecrire les gerbes
qui sont effectivement triviales sur $M$ en terme d'extensions centrales de groupo\"\i des
de $\Gamma\toto M$. Pour d\'ecrire toutes les gerbes sur
$\XX$, on doit passer en g\'en\'eral \`a un groupo\"\i de de Lie Morita-\'equivalent
$\Gamma'\toto M'$.

\par\medskip\centerline{\rule{2cm}{0.2mm}}\medskip
\setcounter{section}{0}

\section{Introduction}

We study $S^1$-bundles and $S^1$-gerbes over differentiable stacks in
terms of Lie groupoids. 

Let $\XX$ be a differentiable stack and $\PP\to\XX$ an $S^1$-bundle
over $\XX$. Let $\Gamma\toto M$ be a Lie groupoid presentation for
$\XX$, i.e., $\XX$ is (isomorphic to) the stack of $\Gamma\toto
M$-torsors. Then $\PP$ gives rise to an $S^1$-bundle $P$ over $M$ on
which $\Gamma\toto M$ acts.  We realize the Chern class of $\PP$ in
terms of connection-like data on $P$ and prove that prequantizations
exist.

Note that $H^2(\Gamma\lcom,\Omega^0)$ contains the obstructions to the
existence of $\PP$ for an arbitrary integer cohomology class and
$H^1(\Gamma\lcom,\Omega^1)$ contains the obstructions to the existence
of a connection on $\PP$ if $\PP$ exists. The possibility of
non-vanishing of these cohomology groups distinguishes our case from
the standard case of manifolds.

If $\GG$ is an $S^1$-gerbe over $\XX$, and $\Gamma\toto M$ a
presentation for $\XX$ as above, then $\GG$ gives rise to a gerbe over
$M$.  So we do not immediately get a description of $\GG$ in terms of
groupoids.  Instead, we can start with a presentation $R\toto M$ of
the differentiable stack $\GG$ and let $\Gamma\toto M$ be the Lie
groupoid defined by the induced presentation $M\to \XX$ of $\XX$, in
other words, $\Gamma=M\times_{\XX}M$.  In this situation, we get a
morphism of groupoids from $R\toto M$ to $\Gamma\toto M$, and,
moreover, $R\to \Gamma$ is an $S^1$-principal bundle.  In fact, $R$ is
an $S^1$-central extension of the Lie groupoid $\Gamma\toto M$. 

Thus the $S^1$-central extensions of $\Gamma\toto M$ are exactly the
$S^1$-gerbes over $\XX$, endowed with a trivialization over $M$.
Therefore, the central extension case is not entirely analogous to
the bundle case. 

Again, we can realize the characteristic class of the gerbe (which we call
the Dixmier-Douady class) in terms of connection-like data and prove
that prequantizations exist.  Note that there are again obstructions
to the existence of honest connective structures and curvings. More
precisely, $H^3(\Gamma\lcom,\Omega^0)$ contains the obstructions to
the existence of $\GG$, given an integer degree-3 cohomology
class. Assuming $\GG$ exists, $H^2(\Gamma\lcom,\Omega^1)$ contains
the obstructions to the existence of a connective structure on
$\GG$. If we assume the existence of a connective structure,
$H^1(\Gamma\lcom,\Omega^2)$ contains the obstructions to the existence
of a curving.

Because groupoid central extensions describe gerbes over $\XX$
together with given trivializations over $M$, we can only describe
those gerbes that are indeed trivial over $M$ in terms of groupoid
central extensions of $\Gamma\toto M$.  To describe all gerbes over
$\XX$, we need to pass in general to a Morita equivalent Lie groupoid
$\Gamma'\toto M'$.

\section{Homology and cohomology}

Let $\Gamma\toto M$ be a Lie groupoid. Define
$\Gamma_p=\underbrace{\Gamma\times_M\ldots\times_M\Gamma}_{\text{$p$
times}}$,
i.e., $\Gamma_p$ is the manifold of composable sequences of $p$ arrows
in the groupoid $\Gamma\toto M$. 
We have $p+1$ canonical maps $\gm_p\to \gm_{p-1}$ (each leaving out
one of the $p+1$ objects involved a sequence of composable arrows),
giving rise to a 
diagram
\begin{equation}\label{sim.ma}
\xymatrix{
\ldots \gm_2
\ar[r]\ar@<1ex>[r]\ar@<-1ex>[r] & \gm_1\ar@<-.5ex>[r]\ar@<.5ex>[r]
&\gm_0\,.}
\end{equation}
In fact,  $\gm \lcom$ is a simplicial manifold. 

The {\em piecewise differentiable chain complex }of $\Gamma\lcom$ is
the total complex 
associated to the double complex $C\lcom(\Gamma\lcom)$. Here
$C_k(\Gamma_p)$ is the free abelian group generated by the
piecewise differentiable maps $\Delta_k\to\Gamma_p$. Its homology
groups $H_k(\Gamma\lcom,\zz)=H_k\big(C\lcom(\Gamma\lcom)\big)$
are called the {\em homology groups }of $\Gamma\toto M$.
 
We denote the dual of the double complex $C\lcom(\Gamma\lcom)$ by
$C\com(\Gamma\lcom)$. It's total cohomology groups 
$H^k(\Gamma\lcom,\zz)=H^k\big(C\com(\Gamma\lcom)\big)$
are called the {\em integer cohomology groups }of $\Gamma\toto M$. 
In the case that $\Gamma\toto M$ is a transformation groupoid $G\times
M\toto M$, these are the $G$-equivariant cohomology groups.  

Finally, we introduce the double complex
$\Omega\com(\Gamma\lcom)$. It's boundary maps are $d:
\Omega^{k}( \gm_p ) \to \Omega^{k+1}( \gm_p )$, the usual exterior
derivative of differentiable forms and $\partial
:\Omega^{k}( \gm_p ) \to \Omega^{k}( \gm_{p+1} )$,  the alternating
sum of the pull back maps of (\ref{sim.ma}).
We denote the total differential by $\delta=(-1)^pd+\del$.
The total cohomology groups of $\Omega\com(\Gamma\lcom)$,
$H_{DR}^k(\Gamma\lcom)=H^k\big(\Omega\com(\Gamma\lcom)\big)$
are called the {\em De~Rham cohomology }groups of
 $\Gamma\toto M$.

Recall that a {\em Morita morphism }from the Lie groupoid
$\Gamma'\toto M'$ to $\Gamma\toto M$ is a morphism of Lie groupoids
satisfying the two conditions
\begin{enumerate}
\item the diagram
$$\xymatrix@=10pt{\Gamma'\dto\rto & M'\times M'\dto\\
\Gamma\rto & M\times M}$$ is cartesian, i.e., a pullback diagram,
\item $M'\to M$ is a surjective submersion.
\end{enumerate}
Two Lie groupoids are Morita equivalent, if and only if there exist a
third Lie groupoid together with a Morita morphism to each of them.

\begin{prop}
Let $f:[\Gamma'\toto M']\to[\Gamma\toto M]$ be a Morita morphism of
Lie groupoids. Then we get induced isomorphisms
$f\upst:H^k(\Gamma\lcom,\zz)\longiso H^k(\Gamma'\lcom,\zz)$
and
$f\upst:H^k_{DR}(\Gamma\lcom)\longiso H^k_{DR}(\Gamma'\lcom)\,.$
\end{prop}

In particular, if $\Gamma\toto M$ is a {\em banal }groupoid, i.e.,
there exists a surjective submersion $\pi:M\to X$, for some manifold
$X$, and $\Gamma\toto M$ is isomorphic to $M\times_XM\toto M$, then we
have canonical isomorphisms
$$f\upst:H^k(X,\zz)\longiso H^k(\Gamma\lcom,\zz)$$
and
\begin{equation}\label{banal}
f\upst:H^k_{DR}(X)\longiso H^k_{DR}(\Gamma\lcom)\,.
\end{equation}

The canonical homomorphism $\Omega\com(\Gamma\lcom)\to
C\com(\Gamma\lcom)\otimes\rr$ induces isomorphisms
\begin{equation}\label{can.iso}
H_{DR}^k(\Gamma\lcom)\longiso H^k(\Gamma\lcom,\rr)\,
\end{equation}
and pairings
$$Z_k(\Gamma\lcom,\zz)\otimes
Z^k_{DR}(\Gamma\lcom)\longrightarrow\rr\,;\quad
\gamma\otimes\omega  \longmapsto {\textstyle\int_{\gamma}\omega}\,.$$

We call a De~Rham cocycle an {\em integer cocycle}, if it
maps under~(\ref{can.iso}) into the image of the canonical map
$H^k(\Gamma\lcom, \zz)\to H^k(\Gamma\lcom,\rr)$. 

\begin{prop}
Let $\omega\in Z^k_{DR}(\Gamma\lcom)$ be a De~Rham cocycle. The
following are equivalent:
\begin{enumerate}
\item  $\omega$ is an integer cocycle,
\item $\int_\gamma\omega\in\zz$, for all $\gamma\in
Z_k(\Gamma\lcom,\zz)$. 
\item for every closed surface $S$ and every $\Gamma\toto M$-torsor
$T$ over $S$, giving rise to a morphism of groupoids $g$ from
$T\times_ST\toto T$ to $\Gamma\toto M$, we have $\int_S
g\upst\omega\in\zz$.  Here we use the isomorphism~(\ref{banal}), to
make sense of the integral.
\end{enumerate}
\end{prop}

(Recall that a {\em $\Gamma\toto M$-torsor }over $S$ is a surjective
submersion $T\to S$, together with an action of $\Gamma\toto M$ on
$T$, such that $S$ is the quotient of $T$ by this action.)

For any abelian sheaf $F$ on the category of differentiable manifolds,
we have the cohomology groups $H^k(\Gamma\lcom,F)$.  One way to define
them is by choosing for every $p$ an injective resolution $F_p\to
I_p\com$ of sheaves on $\Gamma_p$, where $F_p$ is the sheaf induced by
$F$ on $\Gamma_p$; then choosing homomorphisms $f^{-1}I_{p-1}\com\to
I_p\com$ for every map $f:\Gamma_p\to \Gamma_{p-1}$
in~(\ref{sim.ma}). This gives rise to a double complex
$I\com(\Gamma\lcom)$, whose total cohomology groups are the
$H^k(\Gamma\lcom,F)$. 

Examples of abelian sheaves on the category of manifolds are: $\zz$,
$\rr$, $\rr/\zz$,  $\Omega^k$ and $S^1$.  The first three are sheaves
of locally constant functions, $S^1$ is the sheaf of differentiable
$S^1$-valued functions.  With respect to the first three, the notation
$H^k(\Gamma,F)$ does not conflict with the notation introduced
before.

It is well-known that $H^1(\Gamma\lcom,S^1)$ classifies principal
$S^1$-bundles over $\Gamma\lcom$, whereas $H^2(\Gamma\lcom,S^1)$
classifies $S^1$-gerbes over $\Gamma\lcom$.

\section{$S^1$-bundles}

\begin{defn}
Let $\Gamma\toto M$ be a Lie groupoid. A (right) {\em $S^1$-bundle }over
$\Gamma\toto M$ is a (right) $S^1$-bundle $P$ over $M$, together with
a (left) action of $\Gamma$ on $P$, which respects the $S^1$-action,
i.e. we have $(\gamma \cdot x) \cdot t=\gamma \cdot (x \cdot t )$,
for all $t\in S^1$ and  all compatible pairs $(\gamma,x)\in
\Gamma\times_{t,M} P$. 
\end{defn}

Let $Q=\Gamma\times_{t,M} P$ be the manifold of compatible
pairs. Action and
projection form a diagram $Q\toto P$ and it is easy to check that
$Q\toto P$ is in a natural way a groupoid (called the transformation
groupoid of the $\Gamma$-action). Moreover, there is a natural
morphism of groupoids $\pi$ from $Q\toto P$ to $\Gamma\toto M$.  Of
course, $Q$ is an $S^1$-bundle over $\Gamma$.

More is true: the $S^1$-bundle $P$ over $\Gamma\toto M$ gives rise to an
$S^1$-bundle on the simplicial manifold $\Gamma\lcom$.  As such it has
an associated class in $H^1(\Gamma\lcom,S^1)$ and, in fact,
$S^1$-bundles over $\Gamma\toto M$ are classified by
$H^1(\Gamma\lcom,S^1)$.  The exponential 
sequence $\zz\to \Omega^0\to S^1$ induces a boundary map
$H^1(\Gamma\lcom,S^1)\to 
H^2(\Gamma\lcom,\zz)$; the image of the class of $P$ under this
boundary map is called the {\em Chern class }of $P$.

Let $\theta\in \Omega^1(P)$ be a connection form for the $S^1$-bundle
$P\to M$. One checks that $\delta\theta\in C^2_{DR}(Q\lcom)$
descends to $C^2_{DR}(\Gamma\lcom )$. In other words, there exist
unique $\omega\in \Omega^1(\Gamma)$ and $\Omega\in\Omega^2(M)$ such
that $\pi\upst(\omega+\Omega)=\delta\theta$.

\begin{prop}
The class $[\omega+\Omega]\in H^2_{DR}(\Gamma\lcom)$ is
independent of the choice of the connection $\theta$ on $P\to M$. 
Under the canonical homomorphism $H^2(\Gamma\lcom,\zz)\to
H^2_{DR}(\Gamma\lcom)$, the Chern class of $P$ maps to
$[\omega+\Omega]$. 
\end{prop}

Here is a converse.

\begin{prop}
Let $\omega+\Omega\in C^2_{DR} (\Gamma\lcom)$ as above be an integer
2-cocycle.  Then there exists an $S^1$-bundle $P$
over $\Gamma\toto M$ and a connection form $\theta\in \Omega^1(P)$ for
the bundle $P\to M$, such that $\pi\upst(\omega+\Omega)=\delta \theta$.

Moreover, the set of isomorphism classes of all such $(P,\theta)$ is a
simply transitive $H^1 (\Gamma\lcom , \rr/\zz )$-set.  Here
$(P,\theta)$ and $(P',\theta')$ are isomorphic if $P$ and $P'$ are
isomorphic as $S^1$-bundles over $\Gamma\toto M$ and under such an
isomorphism $\theta$ is identified with $\theta'$.
\end{prop}

These two propositions indicate that $\theta$ can be thought of as an
analogue of a connection on $P$ and $\omega+\Omega$ as an analogue of
the curvature of this connection.

On the other hand, we do not call $\theta$ a connection on the
$S^1$-bundle over $\Gamma\toto M$, because this term should be
reserved for $\theta$ satisfying $\del\theta=0$

Thus we suggest the name {\em pseudo-connection }for a connection on
the underlying bundle over $M$.  If $\theta$ is such a
pseudo-connection, we call $\omega+\Omega\in Z^2_{DR}(\Gamma\lcom)$
such that $\pi\upst(\omega+\Omega)=\delta\theta$ the {\em
pseudo-curvature }of $\theta$.

\section{$S^1$-central extensions} 

\begin{defn}
Let ${\gm}\toto {M}$ be a Lie groupoid. An {\em $S^1$-central
extension }of $\gm\toto M$ consists of

1. a  Lie groupoid ${R}\toto {M}$, together with a morphism of Lie
groupoids $(\pi,\id):[R\toto M]\to[\Gamma\toto M]$,

2. a left $S^1$-action on $R$, making $\pi:R\to \Gamma$ a (left)
principal $S^1$-bundle.
\noindent These two structures are compatible in the sense that
$(s\cdot x) (t\cdot y)=st \cdot (xy )$,
for all  $ s,t \in S^{1}$ and $(x, y) \in R\times_MR $.
\end{defn}   

Since $S^1$ is abelian, any left principal $S^1$-bundle is a right
principal $S^1$-bundle in a natural way. Thus, if $R$ and $R'$ are
central extensions of $\Gamma\toto M$ as in the 
definition, we may form the associated bundle  $R\times_{S^1} R'$,
which is again an $S^1$-bundle over $\Gamma$.  It has a natural
groupoid structure making it into another $S^1$-central extension of
$\Gamma\toto M$.  We denote this central extension by $R\otimes R'$.
This operation turns the set of isomorphism classes of $S^1$-central
extensions into an abelian group.

Central extensions of groupoids pull back via morphisms of groupoids.

Groupoid central extensions of $\Gamma\toto M$ give rise to
$S^1$-gerbes over $\Gamma\lcom$, which are trivialized over $M$. Thus
we have the

\begin{prop}
There is a natural exact sequence
\begin{multline*}
H^1(\Gamma\lcom,S^1)\longrightarrow H^1(M,S^1)\longrightarrow \\
\{\text{$S^1$-central extensions of $\Gamma\toto M$}\}
\longrightarrow H^2(\Gamma\lcom,S^1) \longrightarrow H^2(M,S^1)\,.
\end{multline*}
\end{prop}

Given a central extension $R$ of $\Gamma\toto M$, then a connection
form $\theta\in \Omega^1(R)$ for the bundle $R\to \Gamma$, such that
$\del\theta=0$ is a {\em connective structure }on $R$. Given
$(R,\theta)$, a 2-form $B\in\Omega^2(M)$, such that $d\theta=\del B$
is a {\em curving }on $R$, and given $(R,\theta,B)$, the 3-form
$\Omega=dB\in H^0(\Gamma\lcom,\Omega^3)\subset \Omega^3(M)$ is called
the {\em curvature }of $(R,\theta,B)$.  If $\Omega=0$, then
$(R,\theta,B)$ is called a {\em flat } $S^1$-central extension of
$\Gamma\toto M$.  Note that the flat central extensions form an abelian
group.

\begin{prop}
There is a natural exact sequence
\begin{multline*}
H^1(\Gamma\lcom,\rr/\zz)\longrightarrow H^1(M,\rr/\zz)\longrightarrow
\\ \{\text{flat $S^1$-central extensions of $\Gamma\toto M$}\}
\longrightarrow H^2(\Gamma\lcom,\rr/\zz) \longrightarrow
H^2(M,\rr/\zz)\,.
\end{multline*}
\end{prop}

The exponential sequence gives rise to a homomorphism
$H^2(\Gamma\lcom,S^1)\to H^3(\Gamma\lcom,\zz)$.  The image of a
central extension $R$ in $H^3(\Gamma\lcom,\zz)$ is called the {\em
Dixmier-Douady class }of $R$.  The Dixmier-Douady class behaves well
with respect to pullbacks and the tensor operation.

Let $R$ be a central extension of $\Gamma\toto M$. Choose a connection
form $\theta\in\Omega^1(R)$ for the $S^1$-bundle $\pi:R\to
\Gamma$. One checks that $\delta\theta\in Z^3_{DR}(R\lcom)$ descends
to $Z^3_{DR}(\Gamma\lcom)$, i.e., there exist unique
$\eta\in\Omega^1(\Gamma_2)$ and $\omega\in\Omega^2(\Gamma)$ such that
$\pi\upst(\eta+\omega)=\delta\theta$. 

\begin{prop}
The class $[\eta+\omega]\in H^3_{DR}(\Gamma\lcom)$ is independent of
the choice of the connection $\theta$ on $R\to \Gamma$. Under the
canonical homomorphism $H^3(\Gamma,\zz)\to H^3_{DR}(\Gamma\lcom)$, the
Dixmier-Douady class of $R$ maps to $[\eta+\omega]$. 
\end{prop}

Since the class $[\eta+\omega]$ does not change by adding a
coboundary, we may choose, in addition to $\theta$, any
$B\in\Omega^2(M)$, and then the Dixmier-Douady class of $R$ is
represented by $\eta+\omega+\Omega$, such that
$\pi\upst(\eta+\omega+\Omega)=\delta(\theta+B)$. 

\begin{prop}
Given any  3-cocycle  $\eta +\omega+\Omega\in Z^3_{DR}(\Gamma\lcom)$,
as above, satisfying

1. $[\eta+\omega+\Omega]$ is integer,

2. $\Omega$ is exact,

\noindent there exists a groupoid central extension $R\toto M$ of the
groupoid $\gm \toto M$, a connection $\theta$ on the bundle $R\to
\Gamma$ and a 2-form $B\in\Omega^2(M)$, such that
$\delta(\theta+B)=\pi\upst(\eta+\omega+\Omega)$.
The pairs $(R, \theta, B)$ up to isomorphism form a simply transitive 
set under the group of flat central extensions.
\end{prop}

Because of these propositions, $\theta+B$ plays a role similar to a
connection (connective structure plus curving)  on a gerbe over a
manifold. We therefore call $\theta+B$ a {\em pseudo-connection }on
$R$, and $\theta+\omega+\Omega$ its {\em pseudo-curvature}. 

\begin{numrmk}
{\em
Given a 3-cocycle $\eta+\omega+\Omega$ of integer class, we may have
to pass to a Morita equivalent groupoid via a Morita morphism
$[\Gamma'\toto  
M']\to[\Gamma\toto M]$, in order to realize the condition that
$\Omega$ be exact.
For example, if $\Gamma=M$ we may have to pass to an open cover
$\{U_\alpha\}$ of $M$
to construct a groupoid central extension. In this case we use the
Morita morphism $[\coprod_{\alpha,\beta}U_{\alpha\beta}\toto
\coprod_{\alpha}U_\alpha]\to [M\toto M]$. See~\cite{hitchin}. If $M$
is connected, another
possibility is to pass to the (infinite dimensional) path space $PM\to
M$ and use the Morita morphism $[LM\toto PM]\to [M\toto M]$, where
$LM$ is the space of based loops. See~\cite{murray}.}
\end{numrmk}

We close with a theorem that gives an explicit construction of the
central extension with pseudo-connection in the s-connected case.  It
also gives a criterion for a class in $H^3_{DR}(\Gamma\lcom)$ to be
integer. This theorem generalizes the result of~\cite{WX}.

\begin{them}
Let $\gm \toto M$ be an $s$-connected Lie groupoid, and $\eta +\omega
\in C^3_{DR} (\gm\lcom )$ a 3-cocycle, where $\eta \in \Omega^1(\gm_2
)$ and $\omega \in \Omega^2 (\gm )$.  Assume that $\omega $ represents
an integer cohomology class in $H^2_{DR}(\Gamma)$, so that there exists
an $S^1$-bundle $\pi:R\to \gm$ with a connection $\theta \in
\Omega^1 (R)$, whose curvature is $\omega$. 
Assume that $\epsilon \upst R$ endowed with the flat connection
$\epsilon\upst\theta+\pi\upst \epsilon_2\upst\eta$ is holonomy free.  (Here
$\epsilon :M\to\Gamma$ and $\epsilon_2:M\to \Gamma_2$ are the respective identity
morphisms.)
Then $R\toto M$ admits in a natural way the structure of a groupoid,
such that $R$  becomes an $S^1$-central extension of $\Gamma\toto M$
and $\eta +\omega$ the pseudo-curvature of $\theta$.
\end{them}

\Acknowledgements{We thank the Research Institue for the
 Mathematical Sciences in Kyoto,
the Ecole Polytechnique
and Peking University for hospitality and support of the research
summarized in this Note. This research is partially supported by   NSF
          grant DMS00-72171 and NSERC grant 22R81946.       }

%
\end{document}